\documentclass[11pt]{amsart}

\usepackage{amssymb,amsfonts,amsmath}
\usepackage{scalerel}
\usepackage{makecell}
\usepackage{amsthm, stmaryrd}
\usepackage{graphicx}
\usepackage[all,arc]{xy}
\usepackage{hyperref}
\hypersetup{colorlinks,allcolors=violet}
\usepackage{enumerate}
\usepackage{bbm}
\usepackage{dsfont}
\usepackage{mathrsfs}
\usepackage{tikz-cd}
\usetikzlibrary{shapes}
\usepackage{import}
\usepackage{float}
\restylefloat{table}
\usepackage{multirow}
\usepackage{pdflscape}
\usepackage{listofitems}
\usepackage{standalone}
\usepackage[T1]{fontenc}
\usepackage[toc,page]{appendix}

\usepackage[margin=1in]{geometry}
\usepackage{mathtools}
\usetikzlibrary{decorations.pathmorphing}
\usepackage{mathrsfs}
\usepackage{cancel}
\usepackage{relsize}
\newtheorem{thm}{Theorem}[section]
\newtheorem*{thm*}{Theorem}
\newtheorem*{metathm*}{Meta Theorem}

\newtheorem*{setup*}{Setup}

\newtheorem{innercustomthm}{Theorem}
\newenvironment{customthm}[1]
  {\renewcommand\theinnercustomthm{#1}\innercustomthm}
  {\endinnercustomthm}

\newtheorem{innercustomcor}{Corollary}
\newenvironment{customcor}[1]
  {\renewcommand\theinnercustomcor{#1}\innercustomcor}
  {\endinnercustomcor}

  \newtheorem{innercustomprop}{Proposition}
\newenvironment{customprop}[1]
  {\renewcommand\theinnercustomprop{#1}\innercustomprop}
  {\endinnercustomprop}

\newtheorem{cor}[thm]{Corollary}
\newtheorem{prop}[thm]{Proposition}
\newtheorem{lem}[thm]{Lemma}

\theoremstyle{definition}
\newtheorem{defn}[thm]{Definition}

\newtheorem{rem}[thm]{Remark}

\newtheorem*{thm1.2}{\textrm{Theorem 1.2}}

\theoremstyle{remark}
\newcommand{\Mbar}{\overline{\mathcal{M}}}

\newcommand{\M}{\mathcal{M}}

\newtheorem{sch}[thm]{Scholium}

\newcommand{\Z}{\mathbb{Z}}
\newcommand{\Q}{\mathcal{Q}}
\newcommand{\QQ}{\mathbb{Q}}

\renewcommand{\P}{\mathbb{P}}
\newcommand{\bbS}{\mathbb{S}}
\newcommand{\bbL}{\mathbb{L}}

\newcommand{\ch}{\operatorname{ch}}

\newcommand{\cat}[1]{\mathsf{#1}}

\newcommand{\nocontentsline}[3]{}
\let\origcontentsline\addcontentsline
\newcommand\stoptoc{\let\addcontentsline\nocontentsline}
\newcommand\resumetoc{\let\addcontentsline\origcontentsline}

\def\A{\mathbb{A}}

\def\C{\mathbb{C}}

\def\P{\mathbb{P}}

\def\V{\mathbb{V}}

\def\Z{\mathbb{Z}}

\def\calA{\mathcal{A}}

\def\calC{\mathcal{C}}

\def\calM{\mathcal{M}}

\def\calO{\mathcal{O}}
\def\calP{\mathcal{P}}
\def\calQ{\mathcal{Q}}

\def\M{\mathcal{M}}

\newcommand{\Res}{\operatorname{Res}}

\newcommand{\Exp}{\operatorname{Exp}}
\newcommand{\Log}{\operatorname{Log}}

\newcommand{\ve}{\varepsilon}

\newcommand{\DD}{\boldsymbol{D}}
\newcommand{\EE}{\boldsymbol{E}}

\newcommand{\Pic}{\mathrm{Pic}}

\newcommand\cycle[2][\,]{%
  \readlist\thecycle{#2}%
  (\foreachitem\i\in\thecycle{\ifnum\icnt=1\else#1\fi\i})%
}

\newcommand{\bbA}{\mathbb{A}}

\makeatletter
\let\c@equation\c@thm
\makeatother
\numberwithin{equation}{section}

\graphicspath{ {images/} }

\title{Virtual Hodge numbers of $\M_{g, n}(\P^r, d)$: stability and calculations}

\author[S. Kannan]{Siddarth Kannan}\address{Department of Mathematics, Massachusetts Institute of Technology}
\email{\url{spkannan@mit.edu}}
\author[T. Song]{Terry Dekun Song}\address{Department of Pure Mathematics and Mathematical Statistics, University of Cambridge, Cambridge, CB3 0WA}\email{\url{ds2016@cam.ac.uk}}

\begin{document}
\maketitle

\begin{abstract}
We study $\bbS_n$-equivariant motivic invariants of the moduli space $\M_{g, n}(\P^r, d)$ of degree-$d$ maps from $n$-pointed curves of genus $g$ to $\P^r$. In particular, we obtain formulas for the Serre characteristic, which specializes to the Hodge--Deligne polynomial. Fixing $g, r \geq 1$, we prove that an explicit invertible transform of the generating function for the Serre characteristics is rational. We use our formula to prove a stability result for the weight-graded compactly-supported Euler characteristics of $\M_{g, n}(\P^r, d)$ as $d \to \infty$. In genus one and two, we reduce the calculation of the Serre characteristic of $\M_{g, n}(\P^r, d)$ to those of the moduli spaces $\M_{g, n}$ of $n$-pointed curves. Formulas for the latter follow from work of Getzler and Petersen, so our formula in particular determines the Serre characteristic of $\M_{g, n}(\P^r, d)$ for arbitrary $n$, $r$, and $d$ when $g = 1$ and $g = 2$. 
\end{abstract}

\section{Introduction}
Let $\M_{g, n}(\P^r, d)$ denote the moduli space of degree-$d$ maps from smooth $n$-pointed curves of genus $g$ to the projective space $\P^r$. This moduli space admits an action of the symmetric group $\bbS_n$. In this article we study the class $[\M_{g, n}(\P^r, d)]$ in the Grothendieck group of $\bbS_n$-varieties, and then use symmetric function theory to obtain formulas for the $\bbS_n$-equivariant Serre characteristic
\[ \cat{e}^{\bbS_n}(\M_{g, n}(\P^r, d)) = \sum_{i} (-1)^i[H^i_c(\M_{g, n}(\P^r, d);\QQ)] \in K_0(\cat{MHS} ; {\bbS_n}). \]
The $\bbS_n$-equivariant Serre characteristic is the compactly-supported Euler characteristic, considered as an element of the Grothendieck group $K_0(\cat{MHS};\bbS_n)$ of $\bbS_n$-representations in the category $\cat{MHS}$ of mixed Hodge structures over $\QQ$. In particular, the $\bbS_n$-equivariant Serre characteristic determines the Hodge--Deligne polynomial
\begin{equation}\label{eqn:hodgedeligne}
\sum_{i,p,q}(-1)^{i} \mathrm{dim}_\C\left(\mathrm{Gr}^F_{p} \mathrm{Gr}^W_{p + q}H^i_c(\M_{g, n}(\P^r, d);\C)\right)u^p v^q \in \Z[u, v],
\end{equation}
and its $\bbS_n$-equivariant upgrade, in which the coefficient of each $u^pv^q$ in (\ref{eqn:hodgedeligne}) is replaced by a virtual representation of $\bbS_n$. The weight-graded compactly-supported Euler characteristic appearing as the coefficient of $u^p v^q$ in the Hodge--Deligne polynomial is a virtual Hodge number, in the sense that it recovers the usual Hodge number for smooth and proper varieties.

The goal of this paper is to relate the $\bbS_{n}$-equivariant Serre characteristics of $\M_{g,n}(\P^r,d)$ to those of $\M_{g,n}$ (Theorem \ref{thm:Mgreprat}). We completely reduce to $\M_{g, n}$ when $g = 1$ and $g = 2$ (Theorem \ref{thm:lowgenus}). We prove a virtual homological stability result as $d\to \infty$ (Theorem \ref{thm:dstability}). Our results generalize known calculations and stability phenomena in genus zero: see §\ref{subsec:relwork} for more detail.


Our formulas make essential use of the graded $\QQ$-algebra
\[ \hat{\Lambda} := \QQ[\![p_1, p_2, \ldots]\!] \]
of (degree-completed) symmetric functions over $\QQ$. Here $p_i$ is the $i$th power-sum symmetric function, which has degree $i$. The generating function
\begin{equation}\label{eqn:Mgr} \cat{M}_{g, r}:= \sum_{n, d}\cat{e}^{\bbS_n}(\M_{g, n}(\P^r, d))q^d
\end{equation}
is naturally an element of $K_0(\cat{MHS}) \otimes \hat{\Lambda}[\![ q]\!]$. When $g = 1$, the sum in (\ref{eqn:Mgr}) is taken over all $(n, d) \neq (0, 0)$; when $g \geq 2$, the sum is over all $n, d \geq 0$. We will assume that $g \geq 1$, since the genus-zero case is well-understood: see \S\ref{subsec:relwork}. Knowledge of the series $\cat{M}_{g, r}$ is equivalent to knowledge of
\begin{equation}\label{eqn:Mgreps}
    \cat{M}_{g, r}^\varepsilon := \sum_{n, d} \cat{e}^{\bbS_n}(\M_{g, \varepsilon^n}(\P^r, d))q^d,
\end{equation}
where $\M_{g, \varepsilon^n}(\P^r, d)$ is the moduli stack parameterizing degree-$d$ maps to $\P^r$ from $n$-pointed curves, where the $n$ points are no longer required to be distinct. The sum in (\ref{eqn:Mgreps}) is over the same pairs $(n, d)$ as in (\ref{eqn:Mgr}). The two series are related via the plethystic logarithm of $p_1$:
\begin{equation}\label{eqn:logarithm}
\cat{M}_{g, r} = \cat{M}_{g, r}^\varepsilon \circ \Log(p_1),
\end{equation}
where $\circ$ denotes plethysm of symmetric functions. This identity is proved when $r = 0$ in \cite{KSY2024}, and the same proof applies for arbitrary $r$. We prove the rationality in $q$ of a certain invertible transformation of $\cat{M}_{g, r}^{\varepsilon}$. Define the linear operator
\[ \boldsymbol{D}: K_0
(\cat{MHS}) \otimes \hat{\Lambda} [\![ q ]\!] \to K_0
(\cat{MHS}) \otimes \hat{\Lambda} [\![ q ]\!]\]
by
\[ \boldsymbol{D} := \sum_{d > 0} \frac{\partial}{\partial p_d} q^d. \]
We will make use of the operator $\exp(\boldsymbol{D})$, where $\exp$ is the usual exponential. Set $\M_{g, \varepsilon^n} = \M_{g, \varepsilon^n}(\P^0, 0)$, so $\M_{g, \varepsilon^n}$ is the moduli space of smooth genus-$g$ curves together with $n$ ordered, not-necessarily-distinct points. We write
\[ \cat{a}_g^{\varepsilon}:= \sum_{n > 2 - 2g} \cat{e}^{\bbS_n}(\M_{g, \varepsilon^n}) \in K_0(\cat{MHS}) \otimes \hat{\Lambda}. \]
For an integer $d > 0$ define
\[ h_d^\perp : K_0(\cat{MHS}) \otimes \hat{\Lambda} \to K_0(\cat{MHS}) \otimes \hat{\Lambda} \]
by taking $h_d^\perp$ to be the coefficient of $q^d$ in $\exp(\DD)$. This operation is called \textit{skewing} with respect to the homogeneous symmetric function $h_d$; we discuss skewing in more detail in \S\ref{subsec:skew}. Let $\mathbb{L} = [H^2(\P^1;\QQ)] \in K_0(\cat{MHS})$. 
\begin{customthm}{A}\label{thm:Mgreprat}
For all $g \geq 1$, we have
\[ \exp(\boldsymbol{D}) \cat{M}_{g, r}^{\varepsilon} =  f_{g, r}(q) + \frac{h_{2g - 1}^\perp \cat{a}_g^{\ve}}{\bbL^g - 1} \cdot q^{2g - 1} \cdot\frac{\bbL^{r + 1}q - \bbL^{g(r + 1)}q + \bbL^{g(r+1)} - 1}{(1 - q)(1 - \bbL^{r + 1}q)}. \]
where $f_{g, r}(q)$ is a polynomial in $q$ of degree at most $2g - 2$. In particular, $\exp(\DD) \cat{M}_{g, r}^{\varepsilon}$ is a rational function in $q$.
\end{customthm}
Phrased differently, Theorem \ref{thm:Mgreprat} states that $\cat{M}_{g, r}^{\varepsilon}$ is obtained from a rational function in $q$ by applying the invertible transformation $\exp(-\boldsymbol{D})$. We interpret this as a stability result for the Serre characteristic of $\M_{g, n}(\P^r, d)$ as $d \to \infty$. Precisely, we derive the following stability statement for the virtual Hodge numbers of $\M_{g,n}(\P^r,d).$
\begin{customthm}{B}\label{thm:dstability}
Suppose $g,r \geq 1$. Set $h^{p, q}_{vir}(\M_{g, n}(\P^r, d))$ for the coefficient of $u^p v^q$ in the Hodge--Deligne polynomial (\ref{eqn:hodgedeligne}). Then, if
\[ p + q >  2 \cdot ((g + 1)(r + 3) + n + d - 4),\]
we have
\[ h^{p, q}_{vir}\left( \M_{g, n}(\P^r, d + 1) \right) = h^{p - r - 1, q-r - 1}_{vir}\left( \M_{g, n}(\P^r, d) \right). \]
The same statement holds on the level of virtual $\bbS_n$-representations.
\end{customthm}

When $g$, $r$, and $n$ are fixed and $d$ is very large, Theorem \ref{thm:dstability} implies that a proportion of roughly $r/(r + 1)$ of the virtual Hodge numbers of $\M_{g, n}(\P^r, d)$ are stable, in the sense that they are determined by those of $\M_{g, n}(\P^r, d-1).$ This comes from comparing the complex dimension of $\M_{g, n}(\P^r, d),$ which grows linearly in $d$ with slope $r + 1$, and the stable range of weights from the bound given in Theorem \ref{thm:dstability}. The number of stable weights grows linearly in $d$ with slope $2r.$

Following Occam's Razor principle \cite[§1.2]{VakilWood}, it is natural to speculate stabilization of the mixed Hodge structures on the cohomology groups. However, the calculation of the full mixed Hodge structure via Deligne's weight spectral sequence is challenging, for two closely related reasons: the first is the presence of Brill--Noether special line bundles when the mapping degree is small compared to the genus, and the second is the lack of a known smooth normal crossings compactification of $\M_{g,n}(\P^r,d)$ when $g \geq 3$.

In genus one and two, we determine the polynomial $f_{g, r}(q)$ in Theorem \ref{thm:Mgreprat} explicitly in terms of $\cat{a}_{g}^{\ve}$, thereby reducing the calculation of $\cat{M}_{g, r}$ to that of $\cat{a}_g^{\ve}$. We recall the Serre characteristic of $\P^r:$ \[\cat{e}(\P^r) = \bbL^{r} + \bbL^{r - 1} + \cdots + 1. \]

\begin{customthm}{C}\label{thm:lowgenus}
    Let $f_{g, r}(q)$ be as in Theorem \ref{thm:Mgreprat}. When $g = 1$, we have
    \[f_{1, r}(q) = \cat{e}(\P^r) \cdot \cat{a}_{1}^{\ve}, \]
    and when $g = 2$, we have
    \[f_{2, r}(q) = \cat{e}(\P^r)(\cat{a}_2^{\ve} + h_1^\perp \cat{a}_2^{\ve} \cdot q) + \left(\cat{e}(\P^r)\left( \frac{h_{3}^{\perp} \cat{a}_{2}^\ve}{\cat{e}(\P^{2})} - \cat{a}_2^{\ve} \right) + \cat{e}(\P^{2r + 1})\cat{a}_2^{\ve}\right)\cdot q^2. \]
\end{customthm}


The generating function $\cat{a}_{1, r}^{\varepsilon}$ is known via work of Getzler, and the generating function $\cat{a}_{2, r}^{\ve}$ can be derived from work of Getzler and Petersen on local systems on $\M_{2}$, see Remark \ref{rem:a2calculation}. Theorem \ref{thm:lowgenus} thus determines the Serre characteristic $\cat{e}^{\bbS_{n}}(\M_{g, n}(\P^r, d))$ for all $n$, $r$, and $d$ when $g = 1$ or $g = 2$, by taking the plethysm with $\Log(p_1)$ on both sides of the formula in Theorem \ref{thm:Mgreprat}. The theorem implies a very explicit formula for the Serre characteristic of $\M_{1, \varepsilon^{n}}(\P^r, d)$ in terms of those of $\M_{1, \varepsilon^{j}}$ for $n + 1 \leq j\leq n + d - 1$. This formula again involves skewing operations on symmetric functions, this time with respect to the elementary symmetric functions.

\begin{customcor}{D}\label{cor:genusonerecursion}
For $n \geq 1$, let $e_n \in \hat{\Lambda}$ denote the $n$th elementary symmetric function, and let
\[ e_n^\perp : \hat{\Lambda} \to \hat{\Lambda} \]
denote the skewing operation with respect to $e_n$. Then for $d + n > 0$, we have
\[ \cat{e}^{\bbS_n}(\M_{1, \varepsilon^n}(\P^r, d)) = \cat{e}(\P^r)\sum_{k = 1}^{d - 1}(-1)^{k + 1} e_k^\perp \cat{e}^{\bbS_{n + k}} (\M_{1, \ve^{n + k}}) \cdot \left(\bbL^{(d - k)(r + 1)} - \bbL^{(d - k)(r + 1) - r}\right). \]
In particular, we obtain the recursive formula
\[ \cat{e}^{\bbS_n}(\M_{1, \ve^n}(\P^r, d + 1)) = \bbL^{r+1}\cat{e}^{\bbS_n}(\M_{1, \ve^n}(\P^r, d)) +  (-1)^{d + 1}(\bbL^{r + 1} - \bbL)\cdot\cat{e}(\P^r)\cdot e_d^\perp \cat{e}^{\bbS_{n + d}}(\M_{1, \ve^{n + d}}).\]
\end{customcor}

The skewing operators appearing in our results have a concrete representation-theoretic interpretation: for any variety $X$ with an $\bbS_{n + k}$-action,
\[h_k^\perp\cat{e}^{\bbS_{n + k}}(X),\,  e_k^\perp \cat{e}^{\bbS_{n + k}}(X) \in K_0(\cat{MHS}) \otimes \hat{\Lambda} \]
are equal to the multiplicity (considered as a virtual $\bbS_n$-representation) of the  trivial and sign representations of $\bbS_k$ in the $\bbS_{n} \times \bbS_k$-equivariant Serre characteristic of $X$, respectively.



\subsection{Sample calculations}
We include the Serre characteristics of $\M_{1, 1}(\P^2, d)$ and $\M_{2, 0}(\P^2, d)$ for small values of $d$ in Tables \ref{table:m11} and \ref{table:m20}, respectively. We also include the Serre characteristics of the relative Picard stack $\Pic^0_g$ for $g=2$ and $g=3$ in Table \ref{table:pic} (the Serre characteristic of $\Pic^0_g$ is derived as an intermediate result, see \S\ref{subsec:qmappic} and Proposition \ref{prop:picintro} below). In these examples, the Hodge--Deligne polynomial is determined by the substitution $\bbL= uv$. 

In Tables \ref{table:m11} and \ref{table:m20}, we can see the vanishing of the topological Euler characteristic, by taking $\bbL=1$. This vanishing follows from the existence of a fixed-point-free $\C^\star$-action on $\M_{g, n}(\P^r, d)$ for any $g \geq 1$ and $n \geq 0$. By taking $\bbL = 0$, we see that the weight-zero compactly-supported Euler characteristic vanishes in both tables. When $g = 1$, this vanishing follows from the vanishing of the weight-zero subspace of the compactly-supported cohomology of $\M_{1, n}(\P^r, d)$, established in \cite{vzdualcomplex}. In loc. cit., we conjectured that the same vanishing holds for higher $g,$ and Table \ref{table:m20} provides some modest evidence in the case $g = 2$.
 
Table \ref{table:m11} suggests that the stability bound in Theorem \ref{thm:dstability} is not quite optimal. Indeed, when $g = 1$, Corollary \ref{cor:genusonerecursion} implies that the bound can be improved to $p + q > 2\cdot(2r + d + n + 1)$. However, based on the proof of Theorem \ref{thm:dstability}, we strongly expect that the optimal stability bound is still linear in $d$ with slope 2; see Remark \ref{rem:optimality}.




\begin{table}[h]
\begin{tabular}{|c|c|}
\hline
\multicolumn{1}{|l|}{$d$} & $\mathsf{e}(\M_{1, 1}(\mathbb{P}^2, d))$             \\ \hline
$2$                       & $\mathbb{L}^7+2\mathbb{L}^6+\mathbb{L}^5-\mathbb{L}^4-2\mathbb{L}^3-\mathbb{L}^2$                                 \\ \hline
$3$                       & $\mathbb{L}^{10}+2\mathbb{L}^9+\mathbb{L}^8-2\mathbb{L}^7-3\mathbb{L}^6+2\mathbb{L}^4+\mathbb{L}^3-\mathbb{L}^2-\mathbb{L}$                  \\ \hline
$4$                       & $\mathbb{L}^{13}+2\mathbb{L}^{12}+\mathbb{L}^{11}-2\mathbb{L}^{10}-3\mathbb{L}^9+2\mathbb{L}^7-3\mathbb{L}^5-2\mathbb{L}^4+\mathbb{L}^3+2\mathbb{L}^2+\mathbb{L}$ \\ \hline
$5$ & $\bbL^{16}+2\bbL^{15}+\bbL^{14}-2\bbL^{13}-3\bbL^{12}+2\bbL^{10}-3\bbL^8-2\bbL^7+2\bbL^6+4\bbL^5+2\bbL^4-\bbL^3-2\bbL^2-\bbL$ \\ \hline

\end{tabular}
\caption{Serre characteristics of $\M_{1, 1}(\P^2, d)$}
\label{table:m11}
\end{table}

\begin{table}[h]
\begin{tabular}{|c|c|}
\hline
\multicolumn{1}{|l|}{$d$} & $\mathsf{e}(\M_{2, 0}(\mathbb{P}^2, d))$             \\ \hline
$2$                       & $\bbL^8+\bbL^7-\bbL^5-\bbL^4$                                \\ \hline
$3$                       & $\bbL^{10}+\bbL^9-\bbL^7-\bbL^6-\bbL^5-\bbL^4+\bbL^2+ \bbL$                  \\ \hline
$4$                       & $\bbL^{13}+2\bbL^{12}+2\bbL^{11}-\bbL^{10}-3\bbL^9-3\bbL^8-\bbL^7+2\bbL^6+3\bbL^5+2\bbL^4-\bbL^3-2\bbL^2-\bbL$ \\ \hline

\end{tabular}
\caption{Serre characteristics of $\M_{2, 0}(\P^2, d)$}
\label{table:m20}
\end{table}

\subsection{Quasimaps and the relative Picard stack}\label{subsec:qmappic}
Our calculation involves the geometry of an auxiliary space of quasimaps $\Q_{g,n}(\P^r,d)\to \M_{g,n}$, which parametrizes pointed curves together with line bundles and tuples of sections: see Definition \ref{defn:quasimap}. The moduli space $\M_{g,n}(\P^r,d)$ is the open substack of $\Q_{g,n}(\P^r,d)$ parametrizing tuples of sections which do not have common basepoints. Using the stratification of $\Q_{g,n}(\P^r,d)$ by the number of basepoints and  symmetric function theory as in \S \ref{subsec:sym} below, we reduce the calculation of $\cat{M}_{g,r}$  to the analogous quasimap generating function $\cat{Q}_{g,r}:=\sum_{n,d}\cat{e}^{\bbS_n}(\Q_{g,n}(\P^r,d))q^d.$

We then reduce further to the relative Picard stack $\Pic^d_{g,n}\to \M_{g,n}$ of degree-$d$ line bundles. The forgetful map $\Q_{g,n}(\P^r,d)\to \Pic^d_{g,n}$ is a projective bundle when $d\geq 2g-1$ and is in general stratified by projective bundles of varying dimensions, which are controlled by Brill--Noether theory. After an $\bbS_n$-equivariant version of the well-known $d$-independence of $H^*(\Pic^d_{g,n};\QQ)$ (Lemma \ref{lem:pic_independence}), we reduce Theorem \ref{thm:Mgreprat} to the calculation of the Serre characteristic $\cat{e}^{\bbS_n}(\Pic^{2g-1}_{g,n}).$ For this calculation, we use the relative symmetric power $\mathrm{Sym}^{2g - 1}_{g, n} \to \M_{g, n}$, whose fiber over a curve is the symmetric power of the curve. The key point is that
\[\mathrm{Sym}^{2g-1}_{g,n}\to \Pic^{2g-1}_{g,n}\] 
is a projective bundle with fiber $\P^{g-1}$. The generating function for the Serre characteristics of $\mathrm{Sym}^{2g-1}_{g,n}$ is then calculated in terms of the corresponding generating function for $\M_{g, n}$, again using symmetric function theory.

We derive the formulas in Theorem \ref{thm:lowgenus} and Corollary \ref{cor:genusonerecursion} by computing the Grothendieck ring classes of the low-degree ($d < 2g - 1$) quasimap spaces for $g = 1$ and $g =2$. This is immediate when $g = 1$, while the case $g = 2$ involves the geometry of hyperelliptic curves. Finally, we use the change-of-basis between elementary and homogeneous symmetric functions, the $d$-dependence of the Serre characteristics $\cat{e}^{\bbS_n}(\calQ_{g, n}(\P^r, d))$, and dimension bounds provided by Clifford's theorem in order to prove Theorem \ref{thm:dstability}.

Since it may be of independent interest, we record our formula for the Serre characteristic of the relative Picard stack $\Pic^d_{g,n}$.

\begin{customprop}{E}\label{prop:picintro}
    For any $g\geq 1$ and $d\geq 0,$ \[
        \sum_{n\geq 0}\cat{e}^{\bbS_n}(\Pic^d_{g,n}) = \left(\frac{h_{2g-1}^{\perp}(\mathsf{a}_g^{\varepsilon})}{\cat{e}(\P^{g-1})}\right)\circ \Log(p_1).
    \]

\end{customprop}

The formula in Proposition \ref{prop:picintro} is derived from Proposition \ref{prop:Pic_gen_fun}, which concerns the light-weight analogue $\Pic^d_{g,\varepsilon^n}.$ The Serre characteristics of $\Pic_2^0$ and $\Pic_3^0$ are recorded in Table \ref{table:pic} below.

\begin{table}[h]
\begin{tabular}{|c|c|}
\hline
\multicolumn{1}{|l|}{$g$} & $\mathsf{e}(\Pic^0_g)$             \\ \hline
$2$                       & $\bbL^5 + \bbL^4 + \bbL^3 - 1$                                \\ \hline
$3$                       & $\bbL^9 + 2\bbL^8 + 3 \bbL^7 + \bbL^6 + \bbL^3 + 2\bbL^2$
\\ \hline

\end{tabular}
\caption{Serre characteristics of $\Pic^{0}_g$ for $g =2, 3$}
\label{table:pic}
\end{table}


\subsection{The role of symmetric functions}\label{subsec:sym}
Symmetric functions are fundamental to the present work:
\begin{enumerate}
    \item the geometric operation of taking symmetric group quotients for symmetric powers and the basepoint strata in $\Q_{g,n}(\P^r,d)$ corresponds to applying the skewing operator $h_k^{\perp}$ (§\ref{subsec:skew}) on Serre characteristics,

    \item the plethystic logarithm $\Log(p_1)$ appearing in (\ref{eqn:logarithm}) passes between $\mathcal{M}_{g,n}$ and the $n$th fibered power of the universal curve $\mathcal{C}_g \to \M_g$, which is identified with the moduli space of smooth curves with light-weight markings $\mathcal{M}_{g, \varepsilon^n}$ where $\varepsilon<1/n$. Thus markings on the base and the fiber of the universal curve are treated on an equal footing. See \cite[Proposition 8.6]{GetzlerKapranov} for an explicit formula for $\Log(p_1)$.
\end{enumerate}
In Proposition \ref{prop:picintro}, we see that both operators are applied to determine $\mathsf{e}^{\bbS_n}(\Pic^d_{g,n})$ for all $n \geq 0$. We emphasize that the symmetric function formalism is crucial even if one is interested only in the Serre characteristic of the moduli space  $\M_{g,0}(\P^r,d)$ of maps from unmarked curves: conceptually, this is because one still needs to consider the basepoint stratification of the quasimap space; practically, this is clear by first applying $\exp(-\boldsymbol{D})$ to both sides of the formula in Theorem \ref{thm:Mgreprat}, and then extracting the degree $n = 0$ part.

\subsection{Related work}\label{subsec:relwork}
Up to the correction term $f_{g, r}$, Theorem \ref{thm:Mgreprat} expresses $\cat{M}_{g, r}$ in terms of the Serre characteristics of the moduli spaces of curves $\M_{g, n}$, which are of fundamental interest \cite{LooijengaHodge, GetzlerGenusZero, GetzlerMHM, BergstromTommasi, BergstromFaber, BFP, clp}. In the presence of a target $X$ satisfying various assumptions, the topology of the space of maps from a fixed curve $C$ has been studied intensively, especially from the perspective of homological stability. See for example \cite{Segal,cohen,Guest,Banerjee, AA, das_tosteson}.


For many algebro-geometric applications, it is natural to allow the curve to vary in moduli. In this setting, explicit formulas for motivic invariants have so far been limited to maps from rational curves to homogeneous targets. Getzler and Pandharipande computed the Serre characteristic $\cat{e}^{\bbS_n}(\M_{0, n}(\P^r, d))$ for all $n$, $r$, and $d$ in \cite[Theorem 5.6]{GetzlerPandharipande}; they used this calculation to determine the Serre characteristic of the stable maps compactification $\Mbar_{0, n}(\P^r, d)$. See also \cite{Bagnarol} for parallel calculations for genus-zero maps to Grassmannians as well as the recent work \cite{maptoflag} on genus-zero maps to the complete flag variety.

We highlight that the numerical results of \cite{GetzlerPandharipande} indicate the stability property $$\cat{e}^{\bbS_n}(\M_{0, n}(\P^r, d + 1)) = \bbL^{r+1}\cdot \cat{e}^{\bbS_n}(\M_{0, n}(\P^r, d)),$$ which implies the strongest possible upgrade of Theorem \ref{thm:dstability} when $g = 0$. The cohomological refinement of this property is the statement that the Chow and singular cohomology groups of $\mathcal{M}_{0,n}(\P^r,d)$ are independent of $d.$ This has been proven for rational\footnote{On the other hand, work of Cavalieri--Fulghesu \cite{CavalieriFulghesu} describes torsion in the integral Chow ring of $\M_{0,0}(\P^r,d)$ which depends on $d.$ See \S1.2 of their paper for a discussion of this phenomenon.} Chow groups of $\M_{0,0}(\P^r, d)$ and $\M_{0,n}(\P^r,d)$ by Pandharipande \cite{NonlinearGrassmannian} and Oprea \cite[Prop. 2.3.1]{Opthesis} respectively, and for singular cohomology of $\M_{0,n}(\P^r, d)$ by Farb--Wolfson \cite{FarbWolfson}. Oprea has conjectured the stabilization of rational Chow groups of the moduli space of maps from smooth pointed rational curves to flag varieties when the curve class is sufficiently positive\footnote{The assumption on positivity is found to be necessary and added in the online updated version of the paper, see \url{https://mathweb.ucsd.edu/~doprea/flags.pdf}.} \cite[§3]{OpreaFlags}. This body of literature naturally leads to speculations on possible homological stability properties as $d \to \infty$ for $\M_{g,n}(\P^r,d)$ when the genus is positive, which partially inspire the present work. 

It is worth noting that Tosteson \cite{tosteson} has proved a representation stability result for the Kontsevich moduli space $\Mbar_{g, n}(\P^r, d)$ as the number of marked points $n$ tends to infinity. There has also been considerable interest in homological stability for Hurwitz spaces \cite{EVW, LandesmanLevy}, as well as for spaces of rational curves in projective hypersurfaces \cite{circlemethod} and Fano varieties \cite{Manins}, with applications to number theory.

The formulas for $\cat{M}_{g,r}$ in this work are a first step towards calculating Serre characteristics of $\Mbar_{g, n}(\P^r, d)$. Such calculations would be motivic enhancements of our previous works \cite{genusonechar, ks_graph_enum}, which gave an explicit formula for the generating function for $\bbS_n$-equivariant \textit{topological} Euler characteristics of $\Mbar_{1, n}(\P^r, d)$ and a less explicit graph-sum formula for the $\bbS_n$-equivariant Euler characteristics of $\Mbar_{g,n}(\P^r,d)$. In another direction, we expect that homological refinements of our Serre characteristic formulas will involve graph complex techniques, generalizing those which have been effective for $\M_{g,n}$ \cite{cgp, CFGP, PW1journal}.

\subsection*{Acknowledgements}
We are grateful to Dori Bejleri and Davesh Maulik for edifying conversations. We thank Samir Canning, Megan Chang-Lee, Dragos Oprea, Dhruv Ranganathan, Johannes Schmitt, and Ravi Vakil for thoughtful comments on an early draft of this paper. SK is supported by NSF DMS-2401850. TS is supported by a Cambridge Trust international scholarship.

\section{Maps and quasimaps}\label{sec:Picard} Recall that $\M_{g, \varepsilon^n}$ denotes the moduli space of smooth $n$-pointed curves of genus $g$, such that the $n$ ordered points are not necessarily distinct. It can also be thought of as the $n$th fibered power of the universal curve over $\M_g$. The notation is chosen because $\M_{g, \varepsilon^n}$ the complement of the divisor of singular curves in Hassett's compactification $\Mbar_{g, \varepsilon^n}$ of $\M_{g, n}$ \cite{Hassett}, for any weight datum $(\varepsilon, \ldots, \varepsilon) \in \QQ^n$ with $\varepsilon < 1/n$. 
Let \[\mathrm{Pic}^d_{g,\varepsilon^n} \to \M_{g, \varepsilon^n}\]
denote the relative Picard stack of degree $d$ line bundles, and let
\[ \calC_{g, \varepsilon^n} \to \M_{g, \varepsilon^n} \]
denote the universal curve. There is a $\bbS_n$-equivariant isomorphism
\[ \calC_{g, \varepsilon^n} \cong \M_{g, \varepsilon^{n + 1}}, \]
where $\bbS_n$ acts on $\M_{g, \varepsilon^{n + 1}}$ by permuting the last $n$ marked points.

Let $\mathcal{P}_{g, d, n}$ be the Poincaré line bundle on $\mathrm{Pic}_{g,\varepsilon^n}^d\times_{\mathcal{M}_{g,\varepsilon^n}}\mathcal{C}_{g,\varepsilon^n}$: over a point
\[ (C, p_1, \ldots, p_n, L, p_0) \in \mathrm{Pic}_{g,\varepsilon^n}^d\times_{\mathcal{M}_{g,\varepsilon^n}}\mathcal{C}_{g,\varepsilon^n}, \]
the fiber of $\calP_{g, d, n}$ is that of $L$ at $p_0 \in C$. 


Let $$\pi: \mathrm{Pic}_{g,\varepsilon^n}^d\times_{\mathcal{M}_{g,\varepsilon^n}}\mathcal{C}_{g,\varepsilon^n}\to \mathrm{Pic}_{1,\varepsilon^n}^d$$ be the projection map to the first factor. When $d > 2g - 2$ and $r \geq 0$, Riemann--Roch implies that $\pi_*\mathcal{P}_{ g, d, n}^{\oplus r+1}$ forms a vector bundle of rank $(d -g + 1)\cdot(r+1)$ over $\mathrm{Pic}_{g,n}^d$, whose fiber over a point $(C, p_1, \ldots, p_n, L)$ is $H^0(C, L)^{\oplus r+ 1}$. When $d \leq 2g - 2$, the pushforward $\pi_*\calP^{\oplus r+ 1}_{g, d, n}$ is still a coherent sheaf on $\Pic_{g, n}^d$. 

\begin{defn}\label{defn:quasimap}
      For any $d \geq 0$, let $\calQ_{g,\varepsilon^n}(\mathbb{P}^r,d) \to \mathrm{Pic}_{g,\varepsilon^n}^d$ be the total space of the relative projectivization $\mathbb{P}(\pi_*\mathcal{P}_{g, n, d}^{\oplus r+1}).$
\end{defn}
Note that when $d > 2g - 2$, \[\calQ_{g, \ve^n}(\P^r, d) \to \Pic^d_{g, \ve^n}\]
is a Zariski-locally trivial fibration, with fibers $\P^{(d - g + 1)(r+1)}$.
The space $\calQ_{g, \varepsilon^n}(\P^r, d)$ parameterizes tuples $(C, p_1, \ldots, p_n, L, s_0, \ldots, s_r)$ consisting of an $n$-pointed curve $C$, a line bundle $L$ of degree $d$, and $r + 1$ global sections $s_0, \ldots, s_r \in H^0(C, L)$, up to common rescaling. This is the data of a \textit{quasimap} from $C$ to $\P^r$. The mapping space \[\M_{g, \varepsilon^n}(\P^r, d) \hookrightarrow \calQ_{g, \varepsilon^n}(\P^r, d) \]
is the open subset where the sections $s_i$ do not have any common zeros (base points). 

Let
\[ \M_{g, \varepsilon^n}^{(k)}(\P^r, d) \subset \calQ_{g, \ve^n}(\P^r, d) \]
be the locus where the $r + 1$ sections have exactly $k$ common basepoints. Then
\[ \M_{g, \ve^n}(\P^r, d)  = \M_{g, \ve^n}^{(0)}(\P^r, d). \]
Let us denote the Grothendieck group of $\bbS_n$-varieties by $K_0(\cat{Var}_{\C};\bbS_n)$. In this group, we have an equality
\begin{equation}\label{eqn:qmap_stratification}
\sum_{k = 0}^{d}[\M_{g, \ve^n}^{(k)}(\P^r, d)] = [\calQ_{g, \ve^n}(\P^r, d)].  \end{equation}

\begin{lem}\label{lem:basepoint_strata}
There is an $\bbS_n$-equivariant isomorphism of stacks
   \[\M_{g, \varepsilon^n}^{(k)}(\P^r, d) \cong {\M_{g, \varepsilon^{n + k}}^{}(\P^r, d - k)}/{\bbS_k},\]
   where $\bbS_k$ acts by permuting the last $k$ marked points.
\end{lem}
\begin{proof}
Let
\[f: (C, p_1, \ldots, p_n, \{q_i\}_{i = 1}^{k}) \to \P^r \]
be a morphism of degree $d - k$ from a smooth curve $C$ of genus $g$ to $\P^r$, where the $p_i$ and $q_i$ form a set of $n + k$ points of $C$; the $p_i$'s are ordered, but the $q_i$'s are not. The $q_i$'s are equivalently the data of an effective divisor $B = \sum q_i$ of degree $k$ on $C$. Let $\sigma_B$ be the image of the constant function $1$ under
\[ H^0(C, \calO) \to H^0(C, \calO(B)). \]

From the morphism $f$ we get a line bundle
\[L = f^*\calO(1) \in \Pic^{d-k}(C), \]
together with $r + 1$ global sections $s_0, \ldots, s_{r}$ with no common basepoints. We now consider the line bundle
\[ L \otimes \calO(B) \in \Pic^{d}(C) \]
with sections $s_i \otimes \sigma_B \in H^0(C, L\otimes \calO(B))$. This construction defines an isomorphism of stacks
\[ \calM_{g, \varepsilon^{n + k}}(\P^r, d - k)/S_k \to \M_{g, \varepsilon^n}^{(k)}(\P^r, d), \]
as desired.
\end{proof}

\section{Transformations of symmetric functions}\label{sec:Sym}

Consider the graded $\QQ$-algebra
\[ \hat{\Lambda} = \QQ[\![p_1, p_2, \ldots]\!], \]
where the indeterminate $p_i$ has degree $i$. 
We call this algebra the \textit{ring of symmetric functions over $\QQ$}. A basis for the homogeneous degree $n$ part of $\hat{\Lambda}$ is given by the monomials
\[ p_{\mu} := \prod_{i > 0}p_i^{\mu_i} \]
for each partition $\mu \vdash n$, where $\mu_i$ is the number of parts of $\mu$ of size $i$. The homogeneous degree $n$ part of $\hat{\Lambda}$ is identified via the \textit{Frobenius characteristic} with the space of $\QQ$-valued class functions on $\bbS_n$. The Frobenius characteristic identifies the character of an $\bbS_n$-representation $V$ with the element
\[ \ch_n(V) := \frac{1}{n!} \sum_{\sigma \in \bbS_n} \mathrm{Tr}(\sigma) p_{\lambda(\sigma)} \in\hat{\Lambda} \]
where $\lambda(\sigma) \vdash n$ is the cycle type of $\sigma$.

The $\bbS_n$-equivariant Serre characteristic of a variety with $\bbS_n$-action is naturally an element of $K_0(\cat{MHS}) \otimes \hat{\Lambda}$, see e.g. \cite[\S 5]{GetzlerPandharipande} or \cite[Theorem 4.8]{GetzlerPreprint}.
Recall the notation $$\cat{M}_{g,r}^{\varepsilon}=\sum_{n,d} \cat{e}^{\bbS_n}(\mathcal{M}_{g,\varepsilon^n}(\P^r,d))q^d \in K_0(\cat{MHS}) \otimes \hat{\Lambda} [\![ q] \!]$$
as in (\ref{eqn:Mgreps}); the sum is over all $(n, d) \neq (0,0)$ when $g = 1$, and when $g \geq 2$ the sum is over all $n,d\geq 0$.
We now define 
$$\cat{Q}_{g,r}^{\varepsilon}:=\sum_{n\geq 0}\sum_{d\geq 0} \cat{e}^{\bbS_n}(\mathcal{Q}_{g,\varepsilon^n}(\P^r,d))q^d.$$ We will use (\ref{eqn:qmap_stratification}), Lemma \ref{lem:basepoint_strata}, and basic symmetric function theory to pass between these generating functions. Two references for the symmetric function theory that we use are Macdonald \cite{Macdonald} and Getzler--Kapranov \cite{GetzlerKapranov}.
\subsection{Skewing}\label{subsec:skew}
For any \[f = f(p_1, p_2, \ldots) \in \hat{\Lambda}\] we set
\[f^\perp : \hat{\Lambda} \to \hat{\Lambda} \]
for the linear map
\[f \left(\frac{\partial}{\partial p_1}, 2\frac{\partial}{\partial p_2}, 3\frac{\partial}{\partial p_3}, \cdots\right).\]
Sometimes we refer to this operation as \textit{skewing} by $f$.
On the level of $\bbS_n$-equivariant Serre characteristics, taking $\bbS_k$-quotients is encoded by the operation of skewing by the $k$th homogeneous symmetric function $h_k$. For $n > 0$, $h_n$ is homogeneous of degree $n$, and may be defined by the formula
\begin{equation}\label{eqn:hom_pow}1+\sum_{n \geq 1} h_n = \exp\left(\sum_{n \geq 1} \frac{p_n}{n}\right). \end{equation}
Suppose now that $X$ is a variety with an action of $\bbS_{n + k}$. Then, by a standard upgrade of \cite[Proposition 8.10]{GetzlerKapranov} to graded mixed Hodge structures over $\QQ$, we have
\begin{equation}\label{eqn:skew}\cat{e}^{\bbS_n}(X/\bbS_k) = h_k^\perp \cat{e}^{\bbS_{n + k}}(X). 
\end{equation}
To pass between $\cat{M}_{g, r}^{\varepsilon}$ and $\cat{Q}_{g, r}^{\varepsilon}$, we will need to consider the operator $\EE$ on $\hat{\Lambda}[\![ q]\!]$ defined by
\[\EE := \sum_{k \geq 1} h_k^\perp  q^k. \]

Then Lemma \ref{lem:basepoint_strata}, (\ref{eqn:qmap_stratification}), and (\ref{eqn:skew}) together give the following lemma. 
\begin{lem}\label{lem:Qgr-Mgr}
    \begin{align*}
        \EE\,\cat{M}_{g,r}^{\varepsilon} &=  \cat{Q}_{g, r}^{\varepsilon} - \cat{M}_{g, r}^{\varepsilon}.
    \end{align*}
\end{lem}
\begin{proof}
    By definition, \begin{align*}
    \cat{Q}_{g, r}^{\varepsilon} - \cat{M}_{g, r}^{\varepsilon} & = \sum_{d\geq 1}\sum_{k=1}^d \sum_{n\geq 0} [\M_{g,\varepsilon^n}^{(k)}(\P^r,d)]q^d\\ & = \sum_{d\geq 1}\sum_{k=1}^d \sum_{n \geq 0} \left[ \frac{\Res^{S_{n+k}}_{S_n\times S_k}\M_{g, \varepsilon^{n + k}}(\P^r, d-k)}{S_k} \right] q^d\\ & = \sum_{d\geq 1}\sum_{k=1}^d \sum_{n \geq 0} h_k^\perp  \mathsf{e}^{\bbS_{n+k}}(\M_{g, \varepsilon^{n + k}}(\P^r, d-k)) q^d\\ & = \sum_{k\geq 1}\sum_{d'
    \geq 0}\sum_{n'\geq 0}h_k^\perp q^k \cat{e}^{\bbS_{n'}}(\M_{g,\varepsilon^{n'}}(\P^r, d'))q^{d'}\\ & = \EE\,\cat{M}_{g,r}^{\varepsilon}.
\end{align*}
We used change of variables $d':=d-k,n':=n+k$ in the above.
\end{proof}

We now derive an explicit formula for the operator $\boldsymbol{E}$. Recall from the introduction the operator
\[ \boldsymbol{D} = \sum_{n \geq 1} \frac{\partial}{\partial p_n} q^n. \]
\begin{lem}\label{lem:Eformula}
\[ \EE = \exp(\DD) -1. \]
\end{lem}
\begin{proof}
Replace $p_n/n$ by $ \partial/\partial p_n$ in (\ref{eqn:hom_pow}).
\end{proof}
As such, the operator $\exp(\DD)$ passes between $\cat{M}_{g, r}^{\varepsilon}$ and $\cat{Q}_{g, r}^{\varepsilon}$.
\begin{cor}\label{cor:fund_transformation}
    For all $g \geq 1$,
    \[ \exp(\boldsymbol{D}) \cat{M}_{g, r}^{\varepsilon} = \cat{Q}_{g, r}^{\varepsilon}.  \]
    Equivalently,
    \[\cat{M}_{g, r}^{\varepsilon} = \exp(-\boldsymbol{D})\cat{Q}_{g, r}^{\varepsilon}. \]
\end{cor}
 In view of Corollary \ref{cor:fund_transformation}, the statement of Theorem \ref{thm:Mgreprat} is equivalent to the rationality in $q$ of the generating function $\cat{Q}_{g, r}^{\varepsilon}$, which we now pursue.
\begin{defn}
    Define $\cat{Q}_{g,r}^{\varepsilon,+}$ by $$\cat{Q}_{g,r}^{\varepsilon,+}:= \sum_{n\geq 0}\sum_{d\geq 2g-1} \mathsf{e}^{\bbS_n}(\mathcal{Q}_{g,\varepsilon^n}(\P^r,d))q^d.$$ Define $\cat{M}_{g,r}^{\varepsilon, +}$ similarly.
\end{defn}
Of course, rationality of $\cat{Q}_{g, r}^{\varepsilon, +}$ is equivalent to that of $\cat{Q}_{g, r}^{\varepsilon}$. In the next section, we will prove an explicit formula for $\cat{Q}_{g, r}^{\varepsilon, +}$ which in particular is a rational function of $q$.
\section{Symmetric powers and the Picard stack}\label{sec:pic}

The first step is to reduce the calculation of $\cat{Q}_{g, r}^{\varepsilon, +}$ to the Serre characteristic of the relative Picard stack. When $d>2g-2$, we have that \[\mathcal{Q}_{g,\varepsilon^n}(\P^r,d)\to \Pic^{d}_{g,\varepsilon^n}\]
is a projective bundle, with fibers $\P^{(d-g+1)(r+1)-1}$. The Picard stack $\Pic^{d}_{g,\varepsilon^n}$ can then be related to the relative symmetric power \[\mathrm{Sym}^d_{g,\varepsilon^n} \to \M_{g, \varepsilon^n},\] which parametrizes a point in $\M_{g,\varepsilon^n}$ together with $d$ unordered points on the curve. We put
\[ \cat{a}_{g}^{\varepsilon} := \sum_{n \geq 0} \cat{e}^{\bbS_n}(\M_{g, \varepsilon^n}) \in K_0(\cat{MHS}) \otimes \hat{\Lambda} \]
for the generating function for Serre characteristics of $\M_{g, \varepsilon^n}$. 

\begin{defn}
    Define 
    \[
   \exp_{\geq 2g-1}(\DD) : \hat{\Lambda}[\![q]\!] \to  \hat{\Lambda}[\![q]\!] \]
   by
   \[\exp_{\geq 2g-1}(\DD)=\sum_{k\geq 2g-1}h_k^\perp  q^k.\]
\end{defn}
From (\ref{eqn:skew}) we obtain the following identity.
\begin{lem}\label{lem:sym_gen_Fun}
    \[\sum_{n\geq 0}\sum_{d\geq 2g-1}\cat{e}^{\bbS_n}(\mathrm{Sym}^{d}_{g,\varepsilon^n})q^d = \exp_{\geq 2g-1}(\DD) \cat{a}_g^{\varepsilon}.\]
\end{lem}
Our rationality result relies on the following degree independence result of the $\bbS_{n}$-equivariant Serre characteristic of $\Pic_{g, \ve^n}^d$.
\begin{lem}\label{lem:pic_independence}
    The $\bbS_n$-equivariant Serre characteristic $\mathsf{e}^{\bbS_n}(\Pic^{d}_{g,\varepsilon^n})$ is independent of $d$. The same holds for $\mathsf{e}^{\bbS_n}(\Pic^{d}_{g,n}).$
\end{lem}
\begin{proof}
When $g \geq 2$, the argument is exactly the same as \cite[Proposition 0.1]{BMJYJac}. When $g=1$, the argument is easily modified: the map $\Pic^{d}_{1,\varepsilon^n}\to \Pic^{dn}_{1,\varepsilon^n}$ given by $(C, L)\mapsto (C, L^{\otimes n})$ is fiberwise an isogeny over $\M_{1, \ve^n}$, up to translation. By the argument in loc. cit., this induces an $\bbS_n$-equivariant isomorphism
\[H^*(\Pic^{d}_{1,\varepsilon^n}; \QQ) \cong H^*(\Pic^{dn}_{1,\varepsilon^n};\QQ). \]
of mixed Hodge structures. Finally, we apply the $\bbS_n$-equivariant isomorphism \[\Pic^0_{1,\varepsilon^n}\to \Pic^{dn}_{1,\varepsilon^n}\] given by $(C,L)\mapsto (C, L\otimes \mathcal{O}_C(d\sum_{i=1}^n p_i))$. 
\end{proof}

\subsection{Proof of Theorem \ref{thm:Mgreprat}} When $d>2g-2,$ the map $\mathrm{Sym}^{d}_{g,\varepsilon^n}\to \Pic^d_{g, \varepsilon^n}$ taking a degree $d$ divisor $D$ to the line bundle $\calO(D)$ is a Zariski locally trivial projective bundle, with fibers $\P^{d-g}$. We can use this in conjunction with Lemma \ref{lem:basepoint_strata} to derive the following formula for the Serre characteristics of $\Pic^{d}_{g, \varepsilon^n}$, for any $d$.


\begin{prop}\label{prop:Pic_gen_fun}
For any $g \geq 1$ and $d \geq 0$,
\[ \sum_{n \geq 0} \cat{e}^{\bbS_n}(\Pic_{g, \ve^n}^d) = \frac{h_{2g - 1}^\perp \cat{a}_g^{\ve}}{\cat{e}(\P^{g - 1})} \]
and
$$\sum_{n\geq 0}\cat{e}^{\bbS_n}(\Pic^d_{g,\varepsilon^n}) = \left(\exp_{\geq 2g - 1}(\DD) \cat{a}_g^{\varepsilon}\right) \cdot \frac{1 - q}{q^{2g - 1}} \cdot \frac{  (\mathbb{L}-1)(1-\bbL q)}{\bbL q-\bbL^{g}q + \bbL^g -1}.$$
\end{prop}


\begin{proof}
The first formula follows from the fact that $\mathrm{Sym}^{2g-1}_{g,\varepsilon^n}\to \Pic^{2g-1}_{g,\varepsilon^n}$ is a projective bundle with fibers $\P^{g-1}$, together with Lemma \ref{lem:pic_independence}.
    
For the second formula, we use the projective bundle description for all $d\geq 2g-1$, so that \begin{align*}
        \exp_{\geq 2g-1}(\DD) \cat{a}_g^{\varepsilon} & = \sum_{n \geq 0} \sum_{d\geq 2g-1} \cat{e}^{\bbS_n}(\mathrm{Sym}^{d}_{g,\varepsilon^n})q^d\\ & = \sum_{n \geq 0}\cat{e}^{\bbS_n}(\mathrm{Pic}^0_{g,\varepsilon^n})\sum_{d\geq 2g-1}\mathsf{e}(\mathbb{P}^{d-g})q^d.
    \end{align*}
    
    The formula now follows from the calculation
    \begin{align*}
        \sum_{d\geq 2g-1}\mathsf{e}(\mathbb{P}^{d-g})q^d = \frac{q^{2g-1}}{1-q} \cdot \frac{\bbL q-\bbL^{g}q + \bbL^g -1}{(\bbL-1)(1-\bbL q)}.
    \end{align*}
\end{proof}



Equating the two formulas in Proposition \ref{prop:Pic_gen_fun} leads to an identity satisfied by the Serre characteristics of $\M_{g, \ve^n}$, reflecting the relationship between the relative symmetric power and the relative Picard stack.

\begin{cor}\label{cor:combinatorial_identity}
\[ \exp_{\geq 2g - 1}(\DD) \cat{a}_{g}^{\ve} = \left(\frac{h_{2g - 1}^\perp \cat{a}_g^{\ve}}{\cat{e}(\P^{g - 1})} \right)\cdot \frac{q^{2g - 1}}{1 - q} \cdot \frac{\bbL q - \bbL^g q + \bbL^g - 1}{(\bbL - 1)(1 - \bbL q)}. \] In particular, $\exp_{\geq 2g - 1}(\DD) \cat{a}_{g}^{\ve}$ is a rational function in $q.$
\end{cor}

\begin{rem}
Proposition \ref{prop:Pic_gen_fun} should be of independent interest: it reduces the calculation of the Serre characteristic of $\Pic^{d}_{g, \varepsilon^n}$ to that of $\M_{g, \ve^{n + 2g - 1}}$. In particular, we may specialize Proposition \ref{prop:Pic_gen_fun} and apply \cite{CFGP} to obtain the weight-zero compactly-supported $\bbS_n$-equivariant Euler characteristic of $\Pic_{g,\varepsilon^n}$, which we denote by $\chi^{\bbS_n}(W_0H^*_c(\Pic_{g, \ve^n}^d))$. This invariant is equal to the negative of the reduced $\bbS_n$-equivariant topological Euler characteristic of the dual complex of any normal crossings compactification, and the constant term of the Hodge--Deligne polynomial (\ref{eqn:hodgedeligne}). Truncation at weight zero defines a ring homomorphism $K_0(\cat{MHS}) \otimes \hat{\Lambda}\to \hat{\Lambda}$. We thus obtain $$ \chi^{\bbS_n}(W_0 H_c^\star(\Pic^{d}_{g,\varepsilon^n})) = h_{2g - 1}^\perp \chi^{\bbS_{n + 2g - 1}}(W_0H_c^\star(\M_{g,\varepsilon^{n + 2g - 1}})).$$

On the other hand, truncating the formula (\ref{eqn:logarithm}) at weight zero gives $$\sum_{n\geq 0}\chi^{\bbS_n}(W_0H_c^\star(\M_{g,\varepsilon^n})) = \left(\sum_{n\geq 0}\chi^{\bbS_n}(W_0H_c^\star(\M_{g,n}))\right)\circ \Exp(p_1),$$ where $\Exp(p_1) = \sum_{n \geq 1} h_n$ is the plethystic inverse of $\Log(p_1)$ (recall that $h_n$ is defined by (\ref{eqn:hom_pow})). The weight-zero compactly-supported Euler characteristic of $\M_{g,n}$ is computed explicitly in the main theorem of \cite{CFGP}.
\end{rem}




Recall now that $\calQ_{g, n}(\P^r, d)$ is a projective bundle of rank $(d - g + 1)(r + 1) - 1$ over $\Pic_{g, n}^d$ for $d > 2g - 2$. Since a straightforward calculation yields
\[ \sum_{d \geq 2g - 1} \cat{e}(\P^{(d - g + 1)(r + 1)-1}) q^d = \frac{q^{2g - 1}}{1 - q} \cdot \frac{\bbL^{r + 1}q - \bbL^{g(r + 1)}q + \bbL^{g(r + 1)}-1}{(\bbL - 1)(1 - \bbL^{r + 1}q)}, \]
we obtain the following two equivalent formulas via Proposition \ref{prop:Pic_gen_fun}. The second formula is used to prove Theorem \ref{thm:Mgreprat}, and the first formula leads to Corollary \ref{cor:genusonerecursion}.
\begin{prop}\label{prop:rationality}
We have
\[\cat{Q}_{g,r}^{\varepsilon,+} = (\exp_{\geq 2g - 1}(\DD) \cat{a}^{\varepsilon}_g)\frac{(1 - \bbL q)( \bbL^{r + 1}q- \bbL^{g(r + 1)}q   + \bbL^{g(r + 1)}- 1)}{(1 - \bbL^{r+1}q)(\bbL q - \bbL^g q + \bbL^g - 1)}  \]
and
\[ \cat{Q}^{\ve, +}_{g, r} = \frac{h_{2g - 1}^\perp \cat{a}_g^{\ve}}{\bbL^g - 1} \cdot q^{2g - 1} \cdot\frac{\bbL^{r + 1}q - \bbL^{g(r + 1)}q + \bbL^{g(r+1)} - 1}{(1 - q)(1 - \bbL^{r + 1}q)}. \]


In particular, $\cat{Q}_{g,r}^{\varepsilon,+}$ is rational in $q.$
\end{prop}

We have essentially proved Theorem \ref{thm:Mgreprat}. 
\begin{proof}[Proof of Theorem \ref{thm:Mgreprat}]

In view of Corollary \ref{cor:fund_transformation}, Theorem \ref{thm:Mgreprat} is a formula for $\cat{Q}_{g, r}^{\ve}$. This follows from Proposition \ref{prop:rationality}, since the difference $\cat{Q}_{g, r}^\ve -  \cat{Q}_{g, r}^{\ve, +}$ is the polynomial
\[f_{g, r}(q) = \sum_{d = 0}^{2g - 2} q^d\sum_{n \geq 0}\cat{e}^{\bbS_n}(\calQ_{g, \ve^n}(\P^r, d)). \]
Since this polynomial has degree at most $2g - 2$, Theorem \ref{thm:Mgreprat} is proved.
\end{proof}

Expanding the second formula of Proposition \ref{prop:rationality} as a power series in $q$, or by using Proposition \ref{prop:Pic_gen_fun}, we obtain the following explicit formula for the Serre characteristic of $\calQ_{g, n}(\P^r, d)$ when $d > 2g - 2$.

\begin{cor}\label{cor:explicit_stable_qmap}
For $d > 2g - 2$, we have the formula
\[ \cat{e}^{\bbS_n}(\calQ_{g, \ve^n}(\P^r, d)) = h_{2g - 1}^\perp \cat{e}^{\bbS_n}(\M_{g, \ve^{n + 2g - 1}}) \cdot \frac{\bbL^{(r + 1)(d - g + 1)} - 1}{\bbL^{g} - 1}. \]
\end{cor}
Immediate from Corollary \ref{cor:explicit_stable_qmap} is the following recursive formula for the Serre characteristics of quasimap spaces of large degree.
\begin{cor}\label{cor:stable_qmap_diff}
Suppose $d > 2g - 2$. Then
\[\cat{e}^{\bbS_n}(\calQ_{g, \ve^n}(\P^r, d + 1)) - \bbL^{r+1}\cat{e}^{\bbS_n}(\calQ_{g, \ve^n}(\P^r, d)) = h_{2g - 1}^{\perp} \cat{e}^{\bbS_{n + 2g - 1}}(\M_{g, \ve^{n + 2g - 1}}) \cdot \frac{\bbL^{r + 1} - 1}{\bbL^{g}  - 1}. \]
\end{cor}

We now prove Theorem \ref{thm:dstability}, which we restate for the reader's convenience.
\begin{customthm}{B}
Suppose $g,r \geq 1$. Set $h^{p, q}_{vir}(\M_{g, n}(\P^r, d))$ for the coefficient of $u^p v^q$ in the Hodge--Deligne polynomial (\ref{eqn:hodgedeligne}). Then, if
\[ p + q >  2 \cdot ((g + 1)(r + 3) + n + d - 4),\]
we have
\[ h^{p, q}_{vir}\left( \M_{g, n}(\P^r, d + 1) \right) = h^{p - r - 1, q-r - 1}_{vir}\left( \M_{g, n}(\P^r, d) \right). \]
The same statement holds on the level of virtual $\bbS_n$-representations.
\end{customthm}
\begin{proof}
Let $\delta$ denote the degree of the difference between the Hodge--Deligne polynomials of \[{\M_{g, \ve^{n}}(\P^r, d + 1)}\quad \mbox{and} \quad \bbA^{r + 1} \times\M_{g, \ve^{n}}(\P^r, d).\]
We will show that \[\delta \leq 2 \cdot ((g + 1)(r + 3) + n + d - 4),\]
and then use a basic fact about plethysm of symmetric functions to pass from $\M_{g, \ve^n}(\P^r, d)$ to $\M_{g, n}(\P^r, d)$.
Set
\[ A_g(k, d) := (-1)^ke_k^\perp\left(\cat{e}^{\bbS_{n +k}}(\calQ_{g, \ve^{n + k}}(\P^r, d + 1 - k)) - \bbL^{r + 1}\cat{e}^{\bbS_{n + k}}(\calQ_{g, \ve^{n + k}}(\P^r, d - k)) \right). \]
We have the following equality from Corollary \ref{cor:fund_transformation}:
\begin{align*}
\cat{e}^{\bbS_n}(\M_{g, \ve^n}(\P^r, d + 1))& - \bbL^{r + 1}\cat{e}^{\bbS_{n}}(\M_{g, \ve^n}(\P^r, d))\\& = (-1)^{d + 1}\cat{e}(\P^r) e_{d + 1}^\perp \cat{e}^{\bbS_{n + d + 1}}(\M_{g, \ve^{n + d + 1}})
+ A_g(d, d) + \sum_{k = 1}^{d-1} A_g(k, d).
\end{align*}

If $v(X)$ denotes the degree of the Hodge--Deligne polynomial of $X$, we have the bound
\[v(X) \leq 2 \dim(X). \]
We also have
\[v(X) + v(Y) \leq \max(v(X), v(Y)). \]
Let
\[M_{\mathrm{stab}} := \max \left\{v( \calQ_{g, \ve^{n + k}}(\P^r, d - k + 1)) - v(\A^{r + 1}\times\calQ_{g, \ve^{n + k}}(\P^r, d - k) ) \mid 1 \leq k \leq d - 2g + 1 \right\}, \]
\[M_{\mathrm{unstab}} := \max \left\{v( \calQ_{g, \ve^{n + k}}(\P^r, d - k + 1) )- v(\A^{r + 1} \times \calQ_{g, \ve^{n + k}}(\P^r, d - k) ) \mid d - 2g + 2 \leq k \leq d - 1 \right\},\]
\[M_{1}:= v(\P^r \times M_{g, \ve^{n + d + 1}}) , \]
and
\[M_2 :=  \max(v(\calQ_{g, \ve^{n + d}}(\P^r, 1)),  v(\bbA^{r + 1} \times \P^{r} \times \M_{g, \ve^{n + d}})) \]
then
\begin{align*}
\delta \leq \max\{M_{\mathrm{stab}}, M_{\mathrm{unstab}}, M_1, M_2\}.
\end{align*}
We now bound each term separately. \\
\textbf{Bounding $M_{\mathrm{stab}}$:}
In this case we can use Corollary \ref{cor:stable_qmap_diff}, since the bound on $k$ implies that $d - k > 2g - 2$. Since
\[\cat{e}^{\bbS_n}(\Pic^0_{g, \ve^n}) = h_{2g - 1}^{\perp} \cat{e}^{\bbS_{n + 2g - 1}}(\M_{g, \ve^{n + 2g - 1}}) \cdot \frac{\bbL - 1}{\bbL^{g}  - 1},   \]
we can write
\[M_{\mathrm{stab}} \leq \max\{ v\left(\P^r \times \Pic_{g, \ve^{n + k}}^0\right) \mid  1 \leq k \leq d - 2g + 1 \}.  \]
The maximum is clearly attained when $k = d - 2g + 1$. Hence we have
\begin{align*}
    \frac{1}{2}M_{\mathrm{stab}} &\leq r + \dim\left(\Pic_{g, \ve^{n + d - 2g + 1}}^0\right)\\&= r  + 4g - 3 + n + d -2g + 1 \\&= 2g + n + d + r - 2.
\end{align*}
\\
\textbf{Bounding $M_{\mathrm{unstab}}$:}
To bound $M_{\mathrm{unstab}}$, we use Clifford's theorem, which states that for a line bundle $L$ of degree $d$ on a curve $C$, we have $\dim H^0(C, L) \leq \frac{d}{2} + 1$. This implies the dimension bound
\begin{align*}
\dim \left(\bbA^{r + 1} \times \calQ_{g, \ve^{n + k}}(\P^r, d - k)\right) &\leq \left(\frac{d - k}{2} + 2\right)\left( r + 1\right) - 1 + \dim \Pic^{d - k}_{g, \ve^{n + k}} \\&= \left(\frac{d - k}{2} + 2\right)\left( r + 1\right) - 1 + 4g - 3 + n + k.
\end{align*}
which holds at least for $d - 2g + 2\leq k \leq d - 1$. For $k$ in the given range, we see that this bound is maximized when $k = d - 2g + 2$, so we obtain
\begin{align*}
    \dim (\bbA^{r + 1} \times \calQ_{g, \ve^{n + k}}(\P^r, d - k)) &\leq (g + 1)(r + 1) - 1 + 4g - 3 + n + d - 2g + 2\\&= (g + 1)(r + 1) +2g + n + d-2 \\&=(g + 1)(r + 3) + n + d - 4,
\end{align*}
which holds for all $d - 2g + 2\leq k \leq d - 1$. Now we consider $\dim \calQ_{g, \ve^{n + k}}(\P^r, d - k + 1)$ in this range. When $k = d - 2g + 2$, we have
\begin{align*}
\dim \calQ_{g, \ve^{n + d - 2g + 2}}(\P^r, 2g - 1) &= g(r + 1) - 1 + \dim \Pic^{2g - 1}_{\ve^{n + d - 2g + 2}} \\&= g(r + 1) - 1 + 4g - 3 + n + d - 2g + 2\\&= g(r + 3) + n + d - 2.
\end{align*}
For $d - 2g + 3 \leq k \leq d - 1$, we again use Clifford's theorem to find
\begin{align*}
\dim \calQ_{g, \ve^{n + k}}(\P^r, d + 1 - k) &\leq \left(\frac{d + 1 - k}{2} +1\right)(r + 1) - 1 + \dim \Pic^{d + 1 -k}_{g, \ve^{n + k}}\\&= \left(\frac{d + 1 - k}{2} +1\right)(r + 1) + 4g - 3 + n + k,
\end{align*}
which is again maximized for the smallest value of $k$, i.e. $k = d - 2g + 3$. Hence
\begin{align*}
\dim \calQ_{g, \ve^{n + k}}(\P^r, d + 1 - k) &\leq g(r + 1) - 1 + 4g - 3 + n + d - 2g  + 3 \\&= g(r + 3) + n + d- 1,
\end{align*}
for all $d - 2g + 3 \leq k \leq d- 1$. Putting these bounds together, we find
\[\frac{1}{2}M_{\mathrm{unstab}} \leq \max\begin{cases}
    (g + 1)(r + 3) + n + d - 4,
    \\ g(r + 3) + n + d - 2
    \\ g(r + 3) + n + d - 1
\end{cases}, \]
and since $r \geq 1$, we conclude the bound
\[\frac{1}{2}M_{\mathrm{unstab}} \leq (g + 1)(r + 3) + n + d - 4.\]
\textbf{Bounding $M_1$ and $M_2$:}
First note that
\[\frac{1}{2} M_1 \leq r + 3g - 3 + n + d + 1 = 3g + n + d +r - 2. \]
To bound $M_2$, note that
\[ \frac{1}{2}M_2 \leq \max\left(\dim \calQ_{g, \ve^{n + d}}(\P^r, 1),  \dim \left(\bbA^{r + 1} \times \P^r \times \M_{g, \ve^{n + d}} \right)\right)  \]
The Clifford bound reads
\begin{align*}
    \dim \calQ_{g, \ve^{n + d}}(\P^r, 1) &\leq r + 1 - 1 + 4g - 3 + n + d\\&= 4g  +n+d +r- 3,
\end{align*}
where we have used that if $\dim H^0(C, L) \leq (1 + 1/2)$ then $\dim H^0(C, L) \leq 1$. Finally, we have
\begin{align*}
\dim{\M_{g, \ve^{n + d}}} + 2r + 1 &= 3g - 3+ n + d + 2r + 1
\\&= 3g + n + d + 2r - 2.
\end{align*}
Putting all of these bounds together, we conclude that
\[\frac{1}{2}M_2 \leq \max(4g + n +d + r - 3 , 3g + n + d + 2r - 2).\]
\textbf{Completing the bound:} Since the bound on $M_{\mathrm{unstab}}$ is larger than that on $M_{\mathrm{stab}}$, and the bound on $M_2$ is larger than that on $M_1$, we can conclude that
    \[\frac{1}{2}\delta \leq \max\begin{cases}
        (g + 1)(r + 3) + n + d - 4\\ 4g + n + d + r - 3\\ 3g + n + d + 2r - 2
    \end{cases}
    \]
Comparing these terms, we find that for any $g, r \geq 1$
\[ \frac{1}{2}\delta \leq (g + 1)(r + 3) + n + d - 4,  \]
as desired. To finish we observe that by the definition of the plethysm operation on  $\hat{\Lambda} [u, v]$, the following statement holds: if $f \in \hat{\Lambda}[u, v]$ is a polynomial of degree $\leq T$ in $u$ and $v$, then $f \circ g$ is also a  polynomial of degree $\leq T$, for any $g \in \hat{\Lambda}[u, v]$ of degree $0$ in $u$ and $v$. This is because $u$ and $v$ act as constants for plethysm: $u \circ g = u$ and $v \circ g = v$ for any $g$. The theorem now follows from (\ref{eqn:logarithm}).
\end{proof}

\begin{rem}\label{rem:optimality}
    The bound on $M_{\mathrm{unstab}}$ in the proof of Theorem \ref{thm:dstability} seems sub-optimal, since it is obtained rather naively from the Clifford bound. However, every bound in the proof is linear with slope $2$ in $d$, so it does not seem likely that the bound in Theorem \ref{thm:dstability} can be improved asymptotically in $d$.
\end{rem}

\section{Low genus calculations}\label{sec:low_genus}
We now specialize Theorem \ref{thm:Mgreprat} to $g = 1$ and $g = 2$, in order to obtain Theorem \ref{thm:lowgenus} and its corollaries. First, note that the $g = 1$ case of Theorem \ref{thm:lowgenus} is immediate from the identity $\cat{Q}_{1, r}^{\ve, +} = \cat{Q}_{1, r}^{\ve} - \cat{e}(\P^r) \cat{a}_1^{\ve}$. We will now obtain Corollary \ref{cor:genusonerecursion} via the first formula in Proposition \ref{prop:rationality}.

\begin{proof}[Proof of Corollary \ref{cor:genusonerecursion}]
We have
\[ \cat{Q}_{1, r}^{\varepsilon, +} = \cat{Q}_{1, r}^{\ve}-\cat{e}(\P^r)\cat{a}_{1}^{\ve}, \]
and $\exp_{\geq 1}(\DD)= \exp(\DD) - 1$. Hence the first formula in Proposition \ref{prop:rationality} reads
\[\cat{Q}_{1, r}^{\varepsilon}  = \cat{e}(\P^r)[(\exp(\DD) - 1) \cat{a}_1^{\ve}] \frac{1 - \bbL q}{1 - \bbL^{r + 1}q} + \cat{e}(\P^r)\cat{a}_{1}^{\ve}.  \]
Applying $\exp(-\DD)$ to both sides and using the $K_0(\cat{MHS}) \otimes \QQ[\![q]\!]$-linearity of this operator, we get
\begin{equation}\label{eqn:genus_one_special}
\cat{M}_{1, r}^{\varepsilon} = \cat{e}(\P^r) \left( \cat{a}_1^\ve \cdot \frac{1 - \bbL q}{1 - \bbL^{r + 1}q} +\left(\exp(-\DD) \cat{a}_1^\ve\right)\left(1 - \frac{1 - \bbL q}{1 - \bbL^{r + 1}q} \right)\right).
\end{equation}
Corollary \ref{cor:genusonerecursion} is the identification of the bidegree $(n, d)$ parts of the two sides of (\ref{eqn:genus_one_special}), where $n$ denotes symmetric function degree, and $d$ denotes $q$-degree. The key input is the definition of the $n$th elementary symmetric function $e_n \in \hat{\Lambda}$ for $n \geq 0$. The element $e_n$ is homogeneous of degree $n$ and is determined in terms of the power sums by the formula
\[ \sum_{n \geq 0} (-1)^ne_n = \exp\left( -\sum_{k \geq 1}\frac{p_k}{k}  \right). \]
In particular, the substitution $p_k \mapsto k \cdot\partial/\partial p_k$ yields the identity
\[  \sum_{d \geq 0} (-1)^de_d^\perp q^d = \exp(-\DD). \]
Since $e_k$ is homogeneous of degree $k$, the degree of $e_k^\perp f$ for $f \in \Lambda$ is $\deg(f) - k$. The corollary now follows from straightforward algebra.
\end{proof}


\subsection{Maps from genus two curves}
We now consider maps from genus two curves: in this case, we have
\[ \cat{Q}_{2, r}^{\ve} = \cat{Q}_{2, r}^{\ve}(0) + \left(\sum_{d = 1}^{2} \sum_{n\ge 0}\cat{e}^{\bbS_n}(\calQ_{2, n}(\P^r, d))\right) + \cat{Q}_{2, r}^{\ve, +}\]
where $\cat{Q}_{2, r}^{\ve, +}$ is explicitly determined in terms of $\cat{a}_{2}^{\ve}$ in Proposition \ref{prop:rationality}, and $\cat{Q}_{2, r}^{\ve}(0) = \cat{e}(\P^r)\cdot \cat{a}_2^{\ve}$. We now determine the two remaining terms. 

\begin{lem}\label{lem:low_deg_genus2}
The following identities hold in the Grothendieck group $K_0(\cat{Var}_{\C};\bbS_n)$ of $\bbS_n$-varieties:
    \[[\Q_{2,\ve^n}(\P^r,1)] = [\M_{2,\ve^{n + 1}}]\cdot [\P^r],\] and 
    \[[\Q_{2,\ve^n}(\P^r,2)] = ([\Pic^2_{2,\ve^n}]-[\M_{2,\ve^n}])\cdot [\P^r] + [\M_{2,\ve^n}]\cdot [\P^{2r+1}].\]
\end{lem}
\begin{proof}
    For both formulas we consider the map $\Q_{2,n}(\P^r,d)\to \Pic^d_{2,n}$ and stratify the map by the dimensions of the fibers. We use $C\in \M_2$ to denote a genus two curve. A degree one line bundle $L$ on $C$ has $h^0(C, L)\geq 1$ if and only if $L\cong \mathcal{O}_C(p)$ for some $p\in C.$ Therefore, the map $\Q_{2,n}(\P^r,1)\to \Pic^1_{2,n}$ is onto the locus $\{\mathcal{O}_C(p)\mid p\in C\}\subset \Pic^1_{2,n},$ which is the image of the immersion $\mathcal{M}_{2,\ve^n}\to \Pic^1_{2,\ve^n}.$ When restricted to this locus, the forgetful map $\Q_{2,\ve^n}(\P^r,1)\to \Pic^1_{2,\ve^n}$ restricts to a projective bundle with fibers isomorphic to $\P^{r}.$

    Now suppose $d = 2$. A degree 2 line bundle $L$ on a genus two curve $C$ underlies a map to projective space $C\to \P^r$ when $h^0(C, L)\geq 2.$ It is known that the only such line bundle is the canonical bundle $\omega_C,$ and $h^0(C, \omega_C)=2$ for all $C\in \mathcal{M}_2.$ Let $\omega: \mathcal{M}_{2,\ve^n}\to \Pic^2_{2,\ve^n}$ be the section given by the canonical bundle. Because the canonical bundle has constant rank on all genus two curves, the restriction of $\mathcal{Q}_{2,\ve^n}(\P^r,2)\to \Pic^{2}_{2,\ve^n}$ to $\omega(\mathcal{M}_{2,n})$ is a projective bundle, with fibers isomorphic to $\P(H^0(C,\omega_{C})^{\oplus r+1}) \cong \P^{2r + 1}$.
    
    A degree two line bundle $L$ on a genus two curve $C$ has $h^0(C,L)=1$ unless $L\cong \omega_C$. Thus on the complement $\Pic^2_{2,\ve^n}\smallsetminus \omega(\mathcal{M}_{2,\ve^n}),$ the map $\Q_{2,\ve^n}(\P^r,2)\to \Pic^2_{2,\ve^n}$ restricts to a projective bundle with fibers isomorphic to $\P^r$. This proves the identity.
\end{proof}

\begin{rem}
    All genus two curves are hyperelliptic, induced by $C\to |\omega_C|:=\P(H^0(C, \omega_C)).$ From this perspective, we may identify $\P(H^0(C, \omega_C)^{\oplus r+1}) \cong \P(H^0(|\omega_C|, \mathcal{O}(1))^{\oplus r+1})$ and $\mathrm{Map}_2(C, \P^r) \cong \mathrm{Map}_1(\P^1, \P^r).$ Here $\mathrm{Map}_d(\P^1, \P^r)$ is the space of parameterized degree-$d$ maps $\P^1 \to \P^r$. This gives the formula $[\mathcal{M}_{2,n}(\P^r,2)]=[\mathrm{Map}_1(\P^1, \P^r)]\cdot [\mathcal{M}_{2,n}],$ which agrees with Table \ref{table:m20}.
\end{rem}



Applying the Serre characteristic to the second identity in Lemma \ref{lem:low_deg_genus2}, we obtain the second part of Theorem \ref{thm:lowgenus}.

\begin{prop}\label{prop:genus2formula}
Let $f_{g, r} = \cat{Q}_{g, r}^\ve - \cat{Q}_{g, r}^{\ve,+}$. Then
\[f_{2, r} = \cat{e}(\P^r)(\cat{a}_2^{\ve} + p_1^\perp \cat{a}_2^{\ve} \cdot q) + \left(\cat{e}(\P^r)\left( \frac{h_{3}^{\perp} \cat{a}_{2}^\ve}{\cat{e}(\P^{2})} - \cat{a}_2^{\ve} \right) + \cat{e}(\P^{2r + 1})\cat{a}_2^{\ve}\right)\cdot q^2. \]
\end{prop}

\begin{rem}\label{rem:a2calculation}
    We summarize how to calculate $\cat{a}_{2}^{\varepsilon}.$ Let $\pi:\mathcal{C}\to \M_2$ be the universal curve, and let $\V:=\mathrm{R}^1\pi_*\mathbb{Q}$ be the variation of Hodge structures on $\M_2.$ Using the Leray spectral sequence associated to $\M_{2,\varepsilon^n}\to \M_2,$ we have \begin{align*}
        \cat{a}_2^{\varepsilon} & = \sum_{n=0}^{\infty}\mathsf{e}^{\bbS_n}(\M_2, (\QQ\oplus \V\oplus \QQ(-1))^{\otimes n}) = \left(\sum_{a=0}^{\infty}h_a\right)\left(\sum_{c=0}^{\infty}h_c\bbL^{c}\right)\left(\sum_{b=0}^{\infty} \cat{e}^{\bbS_b}(\M_2, \V^{\otimes b})\right)\\ & = \exp\left(\sum_{k \geq 1}\frac{p_k}{k}(1 +\bbL^k) \right)\left(\sum_{\substack{\eta \text{ integer partition}\\\text{with even parts}}}s_{\eta} \right)\left(\sum_{l\geq m\geq 0}s_{l,m}\cdot \mathsf{e}(\M_2, \V_{l,m})\right),
    \end{align*}
where $s_{\eta}$ and $s_{l,m}$ denote the Schur functions associated to the integer partitions $\eta$ and $(l,m),$ and $\V_{l,m}$ is the irreducible $\mathrm{Sp}_4(\QQ)$-representation with highest weight $(l,m).$ The last equality uses the Schur--Weyl duality for $\mathrm{Sp}_4$ \cite{Brauer1937} as recounted in \cite[Theorem 1.2]{Dipper2008}; the underlying symmetric group representations of the irreducible Brauer algebra modules are calculated by \cite[Theorem 4.1]{HW90}, see also \cite[§3]{Cox2009}.

The calculation of $\mathsf{e}(\M_2, \V_{l,m})$ has been outlined in \cite[§4.5]{BergstromFaber}: let $\calA_2$ be the moduli stack of principally polarized abelian surfaces. The Torelli map $t: \M_2\to \calA_2$ induces a stratification $\calA_2 = t(\M_2)\sqcup [\calA_1^2/\bbS_2]$, and the cohomology groups of $\V_{l,m}$ on $[\calA_1^2/\bbS_2]$ and $\calA_2$ have been determined by previous works of Petersen \cite{Petersendellip, PetersenA2} respectively, with the latter proving a conjectural formula of Faber--van der Geer \cite{FabervanderGeer}. We also note that the (numerical) Euler characteristics of $\V_{l,m}$ on $\M_2$ have been calculated by Getzler \cite{Getzler_2002}, and parallel calculations on $\M_3$ have been carried out by Bergström--van der Geer \cite{BvdG}.
\end{rem}

\bibliographystyle{amsalpha}
\bibliography{reference}

\end{document}